\documentclass[12pt, numreferences]{kluwer}
\input{amssymb.sty}
\font\tenfrak=eufm10
\font\sevenfrak=eufm7
\font\fivefrak=eufm5
        \newfam\frakfam \def\frak{\fam\frakfam\tenfrak} 
                \textfont\frakfam=\tenfrak
\font\tenDDl=msbm10  
\font\sevenDDl=msbm7 
\font\fiveDDl=msbm5 
        \newfam\DDlfam \def\DDl{\fam\DDlfam\tenDDl} 
                \textfont\DDlfam=\tenDDl 
\scriptfont\DDlfam=\sevenDDl
        \scriptscriptfont\DDlfam=\fiveDDl
\scriptfont\frakfam=\sevenfrak  
        \scriptscriptfont\frakfam=\fivefrak

\long\def\nodo#1{}
\def\genfd{{\bf k}}
\newcommand{\id}{{\rm id}}
\def\namecite#1{{\sc #1}}
 \def\Fhd#1#2{\smash{\mathop{\hbox to 14mm{\rightarrowfill}}
\limits^{\scriptstyle#1}_{\scriptstyle#2}}}

\def\Fhg#1#2{\smash{\mathop{\hbox to 14mm{\leftarrowfill}}
\limits^{\scriptstyle#1}_{\scriptstyle#2}}}

 \def\fhd#1#2{\smash{\mathop{\hbox to 8mm{\rightarrowfill}}
\limits^{\scriptstyle#1}_{\scriptstyle#2}}}

\def\fhg#1#2{\smash{\mathop{\hbox to 8mm{\leftarrowfill}}
\limits^{\scriptstyle#1}_{\scriptstyle#2}}}

\def\diagram#1{\def\normalbaselines{\baselineskip=0pt
\lineskip=0pt\lineskiplimit=0pt}   \matrix{#1}}

\newtheorem{theorem}{Theorem}

\begin{document}

\begin{opening}
\title{Coherent states for Hopf algebras}
\author{Zoran \surname{\v{S}koda}\email{zskoda@irb.hr}}
\institute{ Theoretical Physics Division,
Institute Rudjer Bo\v{s}kovi\'{c}, P.O.Box 180,
 \\ HR-10002 Zagreb}
\runningtitle{Coherent states for Hopf algebras}
\runningauthor{Zoran \v{S}koda}

\begin{abstract}
Families of Perelomov coherent states are defined
axiomatically in the context of unitary representations
of Hopf algebras. A global geometric
picture involving locally trivial noncommutative fibre
bundles is involved in the construction. 
If, in addition, the Hopf algebra has
a left Haar integral, then a 
formula for noncommutative resolution of identity 
in terms of the family of coherent states holds.
Examples come from quantum groups.
\end{abstract}
\classification{AMS classification}{14A22,16W30,14L30,58B32}

\end{opening}
\newtheorem{defn}{Definition}
\newtheorem{prop}{Proposition}
\newtheorem{lem}{Lemma}

Coherent states were at first defined by 
\namecite{Schr\"{o}dinger} in
quantum optics, and later extended by many people
in various frameworks and generalities; 
often in the context of complex Lie groups and their real forms
~(\cite{Ant-review, Kisil:two-ap}). 
\namecite{Perelomov}~(\cite{Perelomov:1972bd,perelomov}) starts with 
a real Lie group $G$,  
and a unitary irreducible representation 
$T : G \rightarrow {\rm Aut}\,V$ 
on a complex Hilbert space $V$. 
Fix a vector $v_0$ in $V$ so that 
Lie subgroup $H \subset G$ is its projective isotropy 
subgroup (i.e. $h \in H$ iff $hv_0$ equals $v_0$ up to
a constant phase). Hence, there is a unitary character
$\chi : H \rightarrow S^1$ such that $hv_0 = \chi(h) v_0$ for
each $h \in H$. For $G$ compact, the representation $T$
extends to a representation of 
the complexification $G^{\DDl C}$ of $G$.

\vskip .02in
A {\bf family of Perelomov coherent vectors}
in $V$ is a family of vectors $\{C(u), u \in G/H \}$,
such that $C([g])= T(g)v_0$ up to a phase.
Coherent {\it states} are projective classes (rays) of 
coherent {\it vectors},
but in practice one often says ``coherent states'' for both notions.
If $V$ is constructed by the method of geometric quantization,
i.e. as the space of holomorphic sections $\Gamma L$ of the corresponding 
quantization line bundle $L$ over $G/H$, then the coherent vectors
may be defined invariantly in terms of 
that line bundle~(\cite{Rawnsley77}). 
For $G$ a compact form of a semisimple Lie group $G^{\DDl C}$,
the details are in Section~\ref{sec:perelRawn} below. 

\vskip .02in
Hopf algebras appear in physics as symmetries
of noncommutative and quantum spaces~(\cite{KlimykSchmud, 
Mack:1992tg, Majid, ManinQGNG, Varilly}). 
Algebra ${\cal O}(G) = \Gamma {\cal O}_G$ 
of regular functions on 
affine algebraic group $G$ are commutative examples
of Hopf algebras~(\cite{Jantzen:AG}) with coproduct
$\Delta : {\cal O}(G)\to {\cal O}(G)\otimes {\cal O}(G)\cong 
{\cal O}(G\times G)$ given by $
(\Delta f)(g_1,g_2) = f(g_1\cdot g_2)$.
In fact, the category of commutative Hopf algebras 
is antiequivalent to the category of 
affine group schemes~(\cite{Jantzen:AG,Schn:lec}). 
Hence, the noncommutative Hopf algebras 
are thought of as (duals to) noncommutative affine group
schemes~(\cite{Drinf:ICM,Schn:lec}; drawback: $\otimes$ is not
a categorical product of noncommutative rings). Actions of affine group
schemes generalize then to the coactions of Hopf algebras,
which can furthermore be ``structure groups'' of
noncommutative fibre bundles. The total space of
such a bundle is either a single algebra (affine case)
or a more complicated system of algebras or categories
with gluing or localizing mechanism to pass between global
and local description.  
Noncommutative fibre bundles
with coacting Hopf algebras playing the role of a structure group
first appeared in now classical work on smash products and
Hopf-Galois extensions. 

Then, {\sc H-J.~Schneider} introduced in ~\cite{Schn-prinHomSp} 
a crucial descent theorem supporting the geometric torsor intuition 
for faithfully flat Hopf-Galois extensions. 
In a study of noncommutative algebras equipped
with differential calculi, {\sc Majid} and {\sc Brzezi\'{n}ski} 
~(\cite{BrzMaj:quantumgauge}) 
discovered a remarkable condition on differential calculi
which enter the definition of principal bundles in that case.
The coherent states on noncommutative projective
homogeneous spaces, exhibited in the present work, seem to need
a bundle theory extended in a different direction.
To this aim, the present author has extended 
the concepts of Zariski {\it locally trivial} principal fibre 
bundles~(\cite{Skoda:thes, Skoda:ban, Skoda:qbun},
and~\cite{Skoda:loc-coinv1-and2}, Part I)
to the setup where both total and base space are
noncommutative (described locally by noncommutative
algebras) and {\it not necessarily affine}. 

Every complex semisimple Lie group $G^{\DDl C}$ 
is an affine algebraic ${\DDl C}$-group,
and $G^{\DDl C}\to G^{\DDl C}/B$ is an algebraic principal
fibration Zariski locally trivialized in a cover by shifts by
action of Weyl group $W$
of the main Bruhat cell~(\cite{Jantzen:AG}).
Noncommutative analogues of such fibrations, 
derived from quantum matrix groups ${\cal G}_q$, 
are recently exhibited~\cite{Skoda:thes, Skoda:ban}. 
The fibrations trivialize in coaction-compatible 
Ore localizations $S_w^{-1}{\cal G}_q$
labeled by the elements $w$ of the Weyl group $W$.
The trivializations are explicitly computed using
an elaborate {\it Ansatz} involving $q$-$w$-Gauss 
decompositions~(\cite{Skoda:ban}, Theorems 9-12; proofs in~\cite{Skoda:thes}
and~\cite{Skoda:loc-coinv1-and2}, II).

In noncommutative case, it is not appropriate to seek for
individual coherent vectors or rays 
in representation space $V$. 
A {\it family} of coherent vectors $C$ 
should be a section of a noncommutative
bundle $V \otimes L_\chi$ over a 
noncommutative ``coset'' space $X$ 
``parametrizing would-be individual'' coherent states,  
where the fiber 
$V = V_\chi = \Gamma L_\chi = {\rm Ind}^{\cal G}_{\cal B} {\DDl C}_\chi $ 
is an analogue of a holomorphically induced representation space,
$\chi$ is an analogue of a character of the inducing subgroup $B$
and $L_\chi$ is an analogue of the Borel-Weil line bundle.
Our noncommutative coset spaces are patched from charts. 
Local descriptions of $X$ and $C$ in different {\it covers} 
by charts are naturally equivalent. 
Earlier studies of coherent states 
for quantum groups~(\cite{JurSt:coh, Sazdjian:1995yg}) 
used computations in a single local chart. 
One of our goals was to show that states
locally computed in~(\cite{Sazdjian:1995yg}) may be
defined {\it a priori}, regardless coordinate choices. 
The main goal was to find a resolution of unity 
in terms of coherent states of compact quantum groups. 

{\bf Notation for Hopf algebras}:
unit map $\eta$, counit $\epsilon$, multiplication $\mu$,
coproduct $\Delta$, antipode $S$ (do not confuse with $S$ 
and $T$ sometimes used for generic Ore subsets in a ring).
We use {\it \namecite{Sweedler}'s notation}: 
for coproduct $\Delta(a) = \sum a_{(1)}
\otimes a_{(2)}$; for (say right) coactions 
$\rho(v) = \sum v_{(0)} \otimes v_{(1)}$;
and their multiplace extensions~(\cite{Majid}).
Ground ring $\genfd$ is any 
commutative unital (in Sections 1-3, later
$\genfd = {\DDl C}$); the category of left ${\cal E}$-modules
for a $\genfd$-algebra ${\cal E}$ is denoted ${}_{\cal E}{\cal M}$.
The ${\cal B}$-comodule analogue has a superscript (${}^{\cal B}{\cal M}$).
The right-hand versions have a {\it right} sub/super-script instead 
(e.g. ${\cal M}^{\cal B}$), and for bi(co)modules we use combinations.

\section{Prerequisites on Ore localization and covers}

Let ${\cal E}$ be an associative unital ring.
For any multiplicative subset $S \subset {\cal E}$,
define category ${\cal C}_l = {\cal C}_l({\cal E},S)$ as follows.
Objects of ${\cal C}_l$ are pairs $(i,Y)$ where $Y$ is a ring
and $i  : {\cal E} \rightarrow Y$ a unital ring homomorphism, 
such that (i) for each $s \in S$, the image $i(s)$ is invertible;
(ii) the set $\{i(s)^{-1}i(r)\,|\, s \in S,\, r \in {\cal E}\}$ 
is a subring of $Y$; (iii) and $i(r) = 0$ iff $sr = 0$ for some $s \in S$.
A morphism $h : (i,Y)\rightarrow (i',Y')$ is a
ring map $h : Y \rightarrow Y'$ such that $h \circ i = i'$.
The left \namecite{Ore} localization of ${\cal E}$ at $S$
is a universal object $(\iota_{\cal E}, S^{-1}{\cal E})$ 
in ${\cal C}_l({\cal E},S)$. 
It exists iff $S$ is a {\bf left Ore set}~(\cite{Skoda:ban, Skoda:lecs}).

We denote by ${\cal E}-{\rm Mod}$ the Abelian category of 
left ${\cal E}$-modules. 
Every Ore localization $(\iota_{\cal E},S^{-1}{\cal E})$ induces 
an exact localization functor by
\[ Q^* = Q^*_S \, : \,{\cal E}-{\rm Mod}\rightarrow S^{-1} {\cal E}-{\rm Mod}
,\,\,\,\,\,\,M \mapsto S^{-1}{\cal E} \otimes_{\cal E} M.\] 
As an example of a localization functor, 
$Q^*$ has a fully faithful right adjoint $Q_*$, 
which is in Ore case exact, equals the restriction of scalars
and has its own right adjoint $Q^!$. 
The adjunction morphism $\iota : {\rm Id}_{{\cal E}-{\rm Mod}}
\rightarrow Q_* Q^*$ is given by 
$\iota_M = \iota_{\cal E} \otimes_{\cal E} 
{\rm id}_M : {\cal E} \otimes_{\cal E} M \cong M \rightarrow Q_* Q^* (M)$.
Denote the composition $Q_* \circ Q^*=:Q$.
It is an endofunctor in ${\cal E}-{\rm Mod}$.
One often denotes $Q_S(M)$ and $Q^*_S(M)$ by $S^{-1}M$.

\vskip .01in
If $S_1$ and $S_2$ are two left Ore sets, then the set $S_1 S_2$ of all
products $s_1 s_2$, where $s_1 \in S$ and $s_2 \in S$, is not necessarily 
multiplicative, but set $S = S_1 \vee S_2$
multiplicatively generated by the union $S_1 \cup S_2$ is left Ore.
Hence the ``double'' localization $Q_{S_1 \vee S_2}$ is well-defined
and $Q_{S_1 \vee S_2}(R)$ has a canonical structure of a ring.
Two 'consecutive localizations' (compositions of functors)
$Q_{S_1} \circ Q_{S_2}$ and $Q_{S_2} \circ Q_{S_1}$,
are {\it not} rings in general, but they play the role in gluing (see below). 
Canonical natural transformations 
$Q_{S_1} \circ Q_{S_2}\rightarrow Q_{S_1 \vee S_2}$
and $Q_{S_2} \circ Q_{S_1}\rightarrow Q_{S_1 \vee S_2}$ are in 
noncommutative case rarely isomorphisms  
('mutual compatibility' of $Q_{S_2}$ with $Q_{S_1}$), 
and $Q_{S_1} \circ Q_{S_2}\not\cong Q_{S_2} \circ Q_{S_1}$,
so one should not count on this. 
Note also the natural maps 
$Q_{S_1}(\iota_M) : Q_{S_1}(M) \rightarrow Q_{S_1} Q_{S_2} (M)$
and $\iota_{S^{-1}_2 M} : Q_{S_2}(M) \rightarrow Q_{S_1} Q_{S_2} (M)$.

For geometrical purposes, one considers families of Ore localizations
$\{ (\iota_\lambda, S^{-1}_\lambda {\cal E})\}_{\lambda \in \Lambda}$
with localization functors $Q^*_\lambda$, $Q_\lambda$. 
We abbreviate 
\[\begin{array}{lc}
{\cal E}_\lambda := S^{-1}_\lambda {\cal E},&
{\cal M}_{\mu\nu\ldots\lambda} := 
S^{-1}_\lambda\ldots S^{-1}_\nu S^{-1}_\mu M,\\
\iota^\mu_{\mu\nu}:= \iota_{{\cal E}_\mu} : {\cal E}_\mu 
\rightarrow {\cal E}_{\mu\nu},&
\iota^\mu_{\mu\nu M} := \iota^\mu_{\mu\nu} \otimes_{\cal E} \id_M,
\\
\iota^\nu_{\mu\nu} 
:= Q_\nu(\iota_\mu) : {\cal E}_\nu \rightarrow {\cal E}_{\mu\nu},&
\iota^\nu_{\mu\nu M} := \iota^\nu_{\mu\nu} \otimes_{\cal E} \id_M.
\end{array}\]

A family of left (right work as well) 
Ore localizations $\{S^{-1}_\lambda {\cal E}\}_{\lambda in \Lambda}$ 
(cf.~\cite{Rosen:88, Skoda:ban, Skoda:qbun}),
{\bf covers} ${\cal E}$ if the fork diagram
\[\diagram{
 {\cal E} \stackrel{\prod_\lambda \iota_\lambda}
\longrightarrow \prod_{\lambda \in \Lambda} S^{-1}_\lambda {\cal E} 
\fhd{\fhd{i_1}{}}{i_2}
\prod_{(\mu,\nu) \in \Lambda \times \Lambda} 
S^{-1}_{\nu}S^{-1}_{\mu} {\cal E} 
}\]
is an equalizer diagram. 
The upper right map $i_1$ is
$\prod_\lambda e^\lambda \mapsto \prod_{\mu\nu} \iota^\mu_{\mu\nu}(e^\mu)$ 
and the lower right map $i_2$ is 
$\prod_\lambda e^\lambda \mapsto \prod_{\mu\nu} \iota^\nu_{\mu\nu}(e^\nu)$.

When this covering condition holds the analogous equalizer 
property extends to other ${\cal E}$-modules (not only $M = {\cal E}$):

{\bf Globalization lemma.} (in wider generality, \cite{Rosen:88} p. 103)
{\it Suppose a finite family of Ore localizations 
$\{S_\lambda^{-1}{\cal E} \}_{\lambda \in \Lambda}$ covers ${\cal E}$.
Then for every left ${\cal E}$-module $M$ the sequence
\[  0 \rightarrow M \rightarrow \prod_{\lambda \in \Lambda} S_\lambda^{-1}M
\rightarrow \prod_{(\mu,\nu) \in \Lambda\times \Lambda} 
S_{\mu}^{-1} S_\nu^{-1} M \]
is exact, where the first morphism is $m \mapsto \prod \iota_{\lambda,M}(m)$
and the second is 
}
\[ \prod_\lambda m_\lambda \mapsto 
\prod_{(\mu,\nu)} (\iota^\mu_{\mu,\nu,M} (m_\mu) - 
\iota^\nu_{\mu,\nu,M} (m_\nu)). \]
Here the order matters: pairs with $\mu = \nu$ may be (trivially) skipped, 
but unlike in the commutative case {\bf we can not 
confine to the pairs of indices with $\mu < \nu$ only}.

\section{Quantum principal bundles using Ore localizations}

Let ${\cal B}$ be a Hopf algebra. An algebra ${\cal E}$ is
a {\bf ${\cal B}$--comodule algebra} 
if it is given with a ${\cal B}$-coaction $\rho$
which is an algebra map~(\cite{KlimykSchmud, Majid, Montg}).
For commutative ${\cal E}$ and ${\cal B}$ this means that the
affine scheme $E = {\rm Spec}\,{\cal E}$ is given a regular action 
of an affine algebraic group $B = {\rm Spec}\,{\cal B}$.
An Ore localization $S^{-1}{\cal E}$ is 
{\bf $\rho$--compatible}~(\cite{Skoda:ban})
if there is a (unique) ${\cal B}$-coaction $\rho_S$ on $S^{-1}{\cal E}$,
making $S^{-1}{\cal E}$ a ${\cal B}$-comodule algebra
such that the localization map 
$\iota_{\cal E}: {\cal E}\rightarrow S^{-1}{\cal E}$ 
is a map of ${\cal B}$-comodule algebras. 
In commutative case, this means
that ${\rm Spec}\,S^{-1}{\cal E}$ is a $B$-invariant Zariski open subscheme
of ${\rm Spec}\,{\cal E}$. {\bf Localized coinvariants} are those
$e$ in $S^{-1}{\cal E}$ for which $\rho_S(e) = e \otimes 1$ 
i.e. the coinvariants for the ``localized'' coaction $\rho_S$.
They form subalgebra 
$(S^{-1}{\cal E})^{{\rm co}{\cal B}} \subset S^{-1}{\cal E}$.

An $({\cal E},{\cal B})$-{\bf Hopf module} 
is an ${\cal E}$-module $M$, with ${\cal B}$-coaction $\rho_M$,
so that $\rho_M(em) = \rho(e)\rho_M(m)$ for all
$e \in {\cal E}$ and $m \in M$. In commutative
case, Hopf modules correspond to $B$-equivariant 
quasicoherent sheaves over ${\rm Spec}\,{\cal E}$.
They form a category commonly denoted by ${}_{\cal E}{\cal M}^{\cal B}$.

A flat localization functor $Q$ on ${\cal E}-{\rm Mod}$
is $\rho$-compatible if there is a (unique) functor
$Q^{{\cal B}}$ on the category ${}_{\cal E}{\cal M}^{\cal B}$
agreeing with $Q$ after forgetting the comodule structures.

\begin{defn} ~(\cite{Skoda:ban, Skoda:qbun}) \label{defn:bundle}
A {\bf Zariski locally trivial principal ${\cal B}$-bundle}
is an ${\cal E}$-comodule algebra $({\cal E},\rho)$
for which there exists a Zariski local trivialization. 
A Zariski {\bf local trivialization} of 
$({\cal E},\rho)$ consists of

$\,\bullet\,\,$ 
a finite {\bf cover} 
$\{(\iota_\lambda, S^{-1}_\lambda{\cal E})\}_{\lambda \in \Lambda}$ 
of ${\cal E}$ by {\bf $\rho$-compatible Ore localizations}, and

$\,\bullet\,\,$
a family 
$\{ \gamma_\lambda : {\cal B} \rightarrow 
S_\lambda^{-1} {\cal E}\}_{\lambda \in \Lambda}$
of ${\cal B}$-comodule algebra maps.
\end{defn}

Here the ${\cal B}$-comodule structure on ${\cal E}_\lambda$ 
is the one induced by $\rho$-compatibility.
Maps $\gamma_\lambda$ are, in commutative
case, induced by trivializing sections (cf. (2) in~\cite{Skoda:ban}),
and we view $\gamma_\lambda$ as an algebraic replacement for
trivializing sections. We discuss some generalizations in~\cite{Skoda:qbun}.

We now sketch how these fibre bundles 
may be understood as indeed being objects over a
quantum quotient space $(X, {\cal O}_X)$.

An additive functor $f^* : {\cal A} \rightarrow {\cal B}$
between Abelian categories is~(\cite{Ros:NcSch})

$\bullet$ {\bf continuous} if it has a right adjoint, say $f_*$;

$\bullet$ {\bf flat} if it is continuous and exact;

$\bullet$ {\bf almost affine} if it is continuous and its right adjoint 
$f_*$ is faithful and exact;

$\bullet$ {\bf affine} if it is almost affine, and its right adjoint
$f_*$ has its own right adjoint, say $f^{!}$.

Consider the category $\underline{CACat}_*$ 
whose objects are 
pairs of the form $({\cal A},{\cal O}_{\cal A})$
where ${\cal A}$ is a (small in a fixed universe) Abelian category, 
and ${\cal O}_{\cal A}$ is an object in ${\cal A}$,
and where ${\rm Hom}({\cal A},{\cal B})$ is the
set of all (additive) continuous 
functors $f^*$ from ${\cal B}$ to ${\cal A}$, 
equipped with a distinguished isomorphism $\phi$
sending $f^*({\cal O}_{\cal B})$ to ${\cal O}_{\cal A}$.
A morphism in $\underline{CACat}_*$ is flat, (almost) affine etc.
if its underlying ``inverse image'' functor is. A cover is
(almost) affine etc. if each morphism in the cover is such.
A {\bf (relative) quasischeme} ${\cal A}$ over ${\cal C}$ is a morphism 
$g : {\cal A}\rightarrow {\cal C}$ in $\underline{CACat}_*$,
for which there is an almost affine cover $\{ Q_\mu^*\}_{\mu \in M}$ 
by flat localizations where all $g_* \circ Q_{\mu*}$ are exact and faithful.
It is a (relative noncommutative) {\bf scheme}~(\cite{Ros:NcSch}) if 
the cover $\{Q_\mu^*\}_{\mu \in M}$ in the definition
can be actually chosen affine, with all $Q^*_\mu\circ g^*$ affine as well. 

If ${\cal E}$ has a trivializing cover 
$\{(\iota_\lambda, S^{-1}_\lambda{\cal E},
\gamma_\lambda)\}_{\lambda \in \Lambda}$, 
which is in the same time the affine cover 
in the definition of a relative scheme over the
(category of modules over) ground ring $\genfd$, then the category
${}_{\cal E}{\cal M}^{\cal B}$ has the structure of
noncommutative scheme $(X,{\cal O}_X)$ over $\genfd$ as well.
Namely, the gluing of {\it Hopf} modules over charts of the cover
reduces to the globalization of ordinary
modules (analogue of the statement
that if we glue equivariant sheaves over charts on a manifold to
a global sheaf then this sheaf is automatically equivariant).
It follows that the localizations $Q^{\cal B}_\lambda$ 
(which exist and are determined by compatibility and $Q_\lambda$) 
form an affine cover of 
${}_{\cal E}{\cal M}^{\cal B}$, and the local triviality
ensures (cf.~\cite{Schn-prinHomSp}, and~\cite{Skoda:ban}, Section 9)
that for each $\lambda$ there is a 
natural \namecite{Schneider}'s equivalence of categories between 
the localized category ${}_{\cal E}{\cal M}^{\cal B}_\lambda$
and the category $({\cal E}_\lambda)^{{\rm co}{\cal B}}-{\rm Mod}$
of modules over the algebra of localized coinvariants.
By descent, its structure sheaf ${\cal O}_{X}$ corresponds to
the family $\{{\cal U}^\lambda\}_{\lambda \in \Lambda}$
of algebras of localized coinvariants 
${\cal U}^\lambda = ({\cal E}_\lambda)^{{\rm co}{\cal B}}$. 
If ${\cal U} := {\cal E}^{{\rm co}{\cal B}}$ 
is the algebra of global coinvariants (``affine quotient''),
then $(X,{\cal O}_X)$ is actually a relative 
scheme over ${}_{\cal U}{\cal M}$, hence {\it a fortiori} over $\genfd$.

\section{Quantum associated bundles}

For any $\genfd$-coalgebra $C$ (e.g. $C = {\cal B}$),
denote by ${\cal M}^{C}$ ($ resp. {}^{C}{\cal M}$)
the category of right (left) ${C}$-{\it co}modules. 
{\bf Cotensor product} is a bifunctor 
$\Box = \Box^{C} : 
{\cal M}^{C} \times {}^{C}{\cal M} \rightarrow {}_{\genfd}{\cal M}$
which is given on objects as 
$N \Box M  := 
{\rm ker} \left(\rho_N \otimes {\rm id}_M - {\rm id}_N \otimes \rho_M\right)$.
The same formula defines the bifunctor
$\Box : {}_{\cal E}{\cal M}^{\cal B} \times {}^{\cal B}{\cal M} \rightarrow 
{}_{\cal E}{\cal M}$. If $D$ is {\it flat} as a $\genfd$-module
(e.g. $\genfd$ is a field),
and $N$ a left $D$- right $C$-bicomodule,
then the cotensor product $N \Box M$ is a $D$-subcomodule of 
$N \otimes_\genfd M$. In particular, if 
$\pi : D \rightarrow C$ is a surjection of
coalgebras then $D$ is a left $D$- right $C$-bicomodule
via $\Delta_D$ and $(\id \otimes \pi) \circ \Delta_D$
respectively, hence ${\rm Ind}^D_C := D \Box^C \_$ 
is a functor from left $C$- to left
$D$-comodules called the {\bf induction} from $C$ to $D$.

Consider for a moment functor 
${\cal E} \otimes_\genfd \_ : {}^{\cal B}{\cal M}
\rightarrow {}_{\genfd}{\cal M}$ (the superscript is intended!).
Given a map $\gamma : {\cal B} \rightarrow {\cal E}$
of ${\cal B}$-comodules for which there is a convolution--inverse 
$\gamma^{-1}$, 
define natural transformations of functors 
$\kappa^\gamma,\bar\kappa^\gamma : 
{\cal E} \otimes_\genfd \_ \rightarrow {\cal E} \otimes_\genfd \_$ by
$\kappa^\gamma_M (\sum_i e_i \otimes m_i) 
= \sum_i e_i \gamma(m_{i(-1)}) \otimes m_{(0)}$
and $\bar\kappa^\gamma = \kappa_{\gamma^{-1}}$. 
If $\gamma$ is a map of ${\cal B}$-comodule algebras, 
then $\gamma^{-1} = \gamma \circ S$. 
In that case, restrict the natural transformation
$\kappa^\gamma$ to the subfunctor 
${\cal E}^{{\rm co}{\cal B}} \otimes_\genfd \_$ and
$\bar\kappa^\gamma$ to the subfunctor ${\cal E} \Box^{\cal B} \_$
and denote the restrictions simply by $\kappa^\gamma |$ and
$\bar\kappa^\gamma |$.
For any natural transformation of 
functors with values $\Phi : F \rightarrow G$ 
in a (say) Abelian category $G$, 
denote by ${\rm Im}\,\Phi : F \rightarrow G$ 
the functor $M \mapsto {\rm Im}\, \Phi_M(F(M))$. 

\begin{lem} 
(a) $\kappa^\gamma \circ \bar\kappa^\gamma = 
 {\rm Id}_{{\cal E}\otimes_\genfd \_}
 = \bar\kappa^\gamma \circ \kappa^\gamma$;

(b) ${\rm Im} (\kappa^\gamma |) = 
{\cal E}\Box^{\cal B} \_$ and ${\rm Im} (\bar\kappa^\gamma |) = 
{\cal E}^{{\rm co}{\cal B}} \otimes_\genfd \_$.
\end{lem}
\begin{pf} (a) follows by calculation: e.g.
the right-hand equality by
\[ 
\begin{array}{lcl}
\bar\kappa^\gamma (\kappa^\gamma (\sum_i f_i \otimes v_i )) & = &
\bar\kappa^\gamma 
\left( \sum_i f_i \gamma (v_{i(-1)}) \otimes v_{i(0)} \right) \\
& = &\sum_i f_i \gamma(v_{i(-2)}) \gamma^{-1} (v_{i(-1)}) \otimes v_{i(0)}\\
& = &\sum_i f_i \epsilon(v_{i(-1)}) \otimes v_{i(0)} 
= \sum_i f_i \otimes v_i .
\end{array} \]
(a) implies that the two equalities in part (b) are equivalent.
By the definition, ${\rm Im} (\kappa^\gamma_M |) $ consists of
the elements of the form $\sum_j f_j \gamma(m_{j(-1)})\otimes
m_{j(0)}$ where $F = f_j \in {\cal E}^{{\rm co}{\cal B}}$
and $m_j \in M$. The coaction axiom, $\rho(f_j) = f_j \otimes 1$,
and the requirement that $\gamma$ is a map of comodules imply 
\[ (\id_{\cal E} \otimes \rho_M) (F) = 
\sum_j f_j \gamma(m_{j(-2)}) \otimes  m_{j(-1)} \otimes m_{j(0)} 
= (\rho_{\cal E} \otimes \id_M) (F)\]
hence $F \in {\cal E}\Box^{\cal B} M$ and
${\rm Im} (\kappa^\gamma_M |)\subset
{\cal E}\Box^{\cal B} M$.

To prove ${\rm Im} (\kappa^\gamma_M |) \supset
{\cal E}\Box^{\cal B} M$ assume the contrary --
$\exists H \in {\cal E}\Box^{\cal B} M\subset {\cal E}\otimes_\genfd M$ 
not in the image. 
Recall that existence of $\gamma$ implies 
(e.g.~\cite{Montg}, Ch.~4; \cite{Majid}; or ~\cite{Skoda:ban}, Sec.~6)
that $\id \otimes \gamma$ followed by multiplication map in ${\cal E}$
is an isomorphism 
${\cal E}^{{\rm co}{\cal B}}\otimes {\cal B}
\stackrel{\cong}{\rightarrow} {\cal E}$
with inverse which we denote here $\phi^\gamma$.
In particular, elements of ${\cal E}$ 
are additive combinations of elements of the form
$f\gamma(\zeta)$ where $f \in {\cal E}^{{\rm co}{\cal B}}$
and $\zeta \in {\cal B}$.
Hence $\exists n \in {\DDl Z}$,
$\exists f_i \in {\cal E}^{{\rm co}{\cal B}}, \zeta_i \in {\cal B},
i = 1,\ldots, n$ such that
$H = \sum_{i} f_i \gamma(\zeta_i) \otimes m_{i}$.

Since the coaction $\rho_{\cal E}$ 
is an algebra map and $\gamma$ an intertwiner, then 
$\rho(f_i \gamma(\zeta_i)) = f_i \gamma(\zeta_{i(1)}) 
\otimes \gamma(\zeta_{i(2)})$. 
Since $H \in {\cal E}\Box^{\cal B} M$ we have 
\begin{equation}\label{eq:phicotH} 
\phi^\gamma(\rho_{\cal E}\otimes \id)(H)
= \phi^\gamma (\id \otimes \rho_M)(H).
\end{equation}
Evaluated and in Sweedler notation equation~\ref{eq:phicotH} reads
\[
\sum_i f_i \otimes \zeta_{i} \otimes m_{i(-1)} \otimes m_{i(0)}
=
\sum_i f_i \otimes \zeta_{i(1)} \otimes \zeta_{i(2)} \otimes m_i.
\]
Apply $\id\otimes\epsilon\otimes\id\otimes\id$ 
to this equation we obtain
\[
\sum_i f_i \epsilon(\zeta_i) \otimes m_{i(-1)}\otimes m_{i(0)}
=\sum_i f_i \otimes \zeta_i \otimes m.
\]
\[
 {\rm Im}\kappa^\gamma_M | \ni
 \sum_i f_i \epsilon(\zeta_i) \gamma(m_{i(-1)})\otimes m_{i(0)}
= \sum_i f_i \gamma(\zeta_i) \otimes m_i
= H. \,\,\,\,\,\Rightarrow\Leftarrow
\]
\end{pf}

By abuse of notation, let 
$\kappa^\gamma | : {\cal E}^{{\rm co}{\cal B}} \otimes_\genfd \_
\rightarrow {\cal E} \Box \_$  
denote also the corestriction of $\kappa^\gamma |$ onto the image
functor ${\cal E}\Box \_$, and alike for $\bar\kappa^\gamma$.
The lemma easily implies
that $\kappa^\gamma |$
is an equivalence of subfunctors
with inverse~$\bar\kappa^\gamma |$. That is,
the pair of natural transformations $(\kappa^\gamma |,
\bar\kappa^\gamma |)$ extends to a pair of
mutually inverse natural autoequivalences of ${\cal E}\otimes_\genfd \_$,
namely $(\kappa^\gamma, \bar\kappa^\gamma)$.
 

We apply this discussion to our localization picture.
For any local trivialization 
$\Lambda = \{ \iota_\lambda, {\cal E}_\lambda , \gamma_\lambda\}_
{\lambda \in \Lambda}$
of ${\cal E}$, we have the natural transformations
$\kappa_\lambda = \kappa^{\gamma_\lambda}, 
\bar\kappa_\lambda = \kappa^{\gamma_\lambda\circ S}$
for all $\lambda \in \Lambda$, and the coproducts
${\cal E}_\lambda \Box M$  
are locally ``identified'' to
${\cal E}^{{\rm co}{\cal B}}_\lambda \otimes_\genfd M$.  
Introduce functor $\Gamma_\lambda \xi_{\_} 
: {}^{\cal B} {\cal M} \rightarrow {}_{\genfd}{\cal M}$
by $M \mapsto \Gamma_\lambda \xi_M = 
{\cal E}^{{\rm co}{\cal B}}_\lambda \otimes_\genfd M$ and similarily,
for consecutive localizations, $\Gamma_{\lambda\mu\ldots}\xi_{\_}$.
Define the natural transformations $\kappa^\lambda_{\lambda\lambda'}$
by the compositions
\[ 
\kappa^\lambda_{\lambda\lambda'M} : \Gamma_\lambda \xi_M
\stackrel{\kappa_{\lambda,M}}\longrightarrow 
{\cal E}_{\lambda} \Box M
\stackrel{\iota^\lambda_{\lambda\lambda'} \Box M }\longrightarrow
{\cal E}_{\lambda\lambda'} \Box M.
\]
and similarily define $\kappa^{\lambda'}_{\lambda\lambda'}$
using $\iota^{\lambda'}_{\lambda\lambda'}$.
Finally, natural transformations
\[ {\cal K}_{\Lambda\Lambda}  
: \prod_{\lambda \in \Lambda} \Gamma_\lambda \xi_{\_}
\rightarrow 
\prod_{(\lambda,\lambda')\in \Lambda \times \Lambda}
{\cal E}_{\lambda\lambda'} \Box \_ 
\]
are defined by ${\rm pr}_{\lambda\lambda'}\circ {\cal K}_{\Lambda\Lambda}
=  \prod_{\lambda \in \Lambda} 
\kappa^{\lambda}_{\lambda\lambda'}\circ {\rm pr}_\lambda$ where 
${\rm pr}_\mu : 
\prod_\lambda \Gamma_\lambda \xi_{\_} \rightarrow \Gamma_\mu \xi_{\_}$
are the natural projection transformations and ${\rm pr}_{\lambda\lambda'}$ 
alike; and similarily define ${\cal K}'_{\Lambda\Lambda}$ 
by using $\kappa^{\lambda'}_{\lambda\lambda'}$. The
{\bf global sections of associated vector bundle functor} 
$\Gamma_\Lambda \xi_{\_}$ is the subfunctor of  
$\prod_{\lambda \in \Lambda} \Gamma_\lambda \xi_{\_}$
such that the fork diagram 
\[
 \diagram{
 \Gamma_\Lambda \xi_{\_} \stackrel{{\rm in}}
\longrightarrow  \prod_{\lambda \in \Lambda} \Gamma_\lambda \xi_{\_}
\fhd{\fhd{{\cal K}_{\Lambda\Lambda}}{}}{{\cal K'}_{\Lambda\Lambda}}
\prod_{(\lambda,\lambda')\in \Lambda \times \Lambda}
{\cal E}_{\lambda\lambda'} \Box \_ 
}
\]
is an equalizer diagram of natural transformations.
\begin{theorem}\label{thm:invarGlobSec}
Functors $\Gamma_\Lambda \xi_{\_}$ are naturally equivalent for
different local trivializations $\Lambda$. There is a natural
equivalence ${\cal K}_\Lambda : \Gamma_\Lambda \xi_{\_}
\rightarrow {\cal E} \Box \_$
making the following diagram sequentially commute:
\begin{equation}\label{eq:qbundlediag}\begin{array}{ccccc}
\Gamma_\Lambda \xi_{\_} & \stackrel{{\rm in}}
\longrightarrow & \prod_{\lambda \in \Lambda} \Gamma_\lambda \xi_{\_}
& \fhd{\fhd{{\cal K}_{\Lambda\Lambda}}{}}{{\cal K'}_{\Lambda\Lambda}}
& \prod_{(\lambda,\lambda')\in \Lambda \times \Lambda}
{\cal E}_{\lambda\lambda'} \Box \_ \\
\,\,\,\,\downarrow {\cal K}_\Lambda
& & \downarrow \prod_\lambda \kappa_\lambda & & 
\,\,\,\,\,\,\,\,\,\,\,\,\,\,\,\| \\
{\cal E} \Box \_ & \stackrel{\iota_\Lambda \Box \_}
\longrightarrow  & \prod_{\lambda \in \Lambda} 
{\cal E}_{\lambda} \Box \_
& \fhd{\fhd{i_1 \Box \_}{}}{i_2 \Box \_}
& \prod_{(\lambda,\lambda')\in \Lambda \times \Lambda}
{\cal E}_{\lambda\lambda'} \Box \_ 
\end{array}
\end{equation}
\end{theorem}
\vskip .1in
\begin{pf}
The square on the right is manifestly commutative by the construction of the
maps involved. Since upper and lower fork diagrams are equalizer diagrams,
and the right vertical arrows natural equivalences,
the transformation ${\cal K}_\Lambda$ exists and 
is uniquely defined by the rest of the diagram.
\end{pf}

As a consequence, given two different local trivializations 
$\Lambda, \Lambda'$ the transformation
${\cal K}_{\Lambda'}\circ {\cal K}_\Lambda
: \Gamma_\Lambda \xi_{\_} \stackrel{\cong}{\longrightarrow}
\Gamma_{\Lambda'} \xi_{\_}$ is a canonical isomorphism of functors.
Hence for fixed $M$, we can denote by $\Gamma \xi_M$
the equivalence class of pairs of the form $(\Lambda,\Gamma_\Lambda\xi_M)$.

Actually, for more general localizations, that is for $Q^*_\mu$,
corresponding to any $\rho$-compatible radical filter ${\cal F}_\mu$, and
given Zariski local trivialization $\Lambda$, we
can define $\Gamma_{\mu\Lambda} \xi_{\_}$ by essentially the same procedure
for ${\cal E}_\mu$ instead of ${\cal E}$, and again, it is 
independent of the choice of $\Lambda$. 
The lattice ${\cal L}_{\cal E}^{\cal B}$
of $\rho$-compatible radical filters is identified
with the lattice ${\frak L}_X$ of flat localizations
of the quotient noncommutative scheme $X$. 
Thus we obtain a bifunctor $\Gamma_{\_} \xi_{\_}$, i.e.
a presheaf  $\mu \mapsto \Gamma_\mu \xi_{\_}$ of functors
$\Gamma_\mu \xi_{\_} : M \mapsto \Gamma_\mu \xi_M$ over ${\frak L}_X$. 
One can consider more general situations where the global
space is not described by only one algebra, but the 
Zariski principal bundles make sense, and the
presheaf  $\mu \mapsto \Gamma_\mu \xi_{\_}$ of functors
over ${\frak L}_X$ still makes sense, whereas the global cotensor
product does not.
Alternatively, we may construct the associated bundles by means
of transition matrices~(\cite{Skoda:qbun}).

{\bf Quantum line bundles.} 
Let us specialize now to $\genfd = {\DDl C}$, choose
a group-like element $\chi\in {\cal B}$
and consider the 1-dimensional left comodule $M = {\DDl C}_\chi$
given by $\rho_M(m) = m \otimes \chi$.
Denote the ``line bundle'' $\xi_M$ by $L_\chi$. 
Its space of sections $\Gamma L_\chi$ can be
identified with a $\DDl C$-subspace of 
$\prod_{\lambda} {\cal E}_\lambda$.
Namely, write explicitly maps ${\cal K}_{\Lambda\Lambda'}$
and  ${\cal K'}_{\Lambda\Lambda'}$, and use the identifications
${\cal E}_\lambda\otimes M \cong {\cal E}_\lambda$ in 
expressing~(\ref{eq:qbundlediag}) to obtain

\begin{equation}\label{eq:gammalchi}
 \Gamma L_\chi \cong \left\lbrace \left.
 f = \prod_{\lambda \in \Lambda} f_\lambda
\,\right\vert\,\begin{array}{c}
f_\lambda \gamma_\lambda(\chi) = f_{\lambda'}\gamma_{\lambda'}(\chi)\\
\,\forall \lambda, \lambda'\,\mbox{ in } {\cal E}_{\lambda\lambda'} 
\mbox{ and in } {\cal E}_{\lambda'\lambda} \end{array}\right\rbrace
\end{equation}


\begin{theorem}\label{thm:globalsec}
$\Gamma L_\chi$ is naturally isomorphic to 
the cotensor product ${\cal E}\Box_{\cal B} M$
as a $\genfd$-vector space.
\end{theorem}

\vskip .1in
Now assume ${\cal E}$ is an algebra 
in category of left $D$- right ${\cal B}$-comodules 
where $D$ is a $\genfd$-coalgebra, flat as $\genfd$-module.
Particular case of importance to us 
is when ${\cal G} = {\cal E} = D$ is a Hopf $\genfd$-algebra
and $\pi : {\cal G} \rightarrow {\cal B}$ is a surjective
homomorphism of Hopf algebras, with natural right ${\cal B}$-comodule
structure $(\id \otimes \pi) \circ \Delta_{\cal G}$. 

Now repeat the arguments preceding Theorem~\ref{thm:invarGlobSec},
but now with the (compatible with other data) 
left $D$-comodule structure added.
More specifically, we consider functor 
${\cal E} \otimes_\genfd \_ : {}^{\cal B}{\cal M}\rightarrow {}^D {\cal M}$,
and the corresponding versions of 
natural transformations $\kappa^\gamma$ etc. and conclude
that the 'new' $\kappa^\gamma$ 
is a natural autoequivalence of the functor
${\cal E} \otimes_\genfd \_ : {}^{\cal B}_{\cal E}{\cal M}
\rightarrow {}^D {\cal M}$ with inverse $\bar\kappa^\gamma$,
inducing the equivalence of subfunctors between
${\cal E}\Box {}^{\cal B}\_$
and ${\cal E}^{{\rm co}{\cal B}}\otimes_\genfd \_$.
Furthermore, construct the left $D$-comodule enrichment of
sections modules $\Gamma_\lambda \xi_{\_}$, $\Gamma_\Lambda \xi_{\_}$
and notice that the morphisms ${\cal K}_{\Lambda\Lambda}$ etc.
in the corresponding commutative diagrams respect the left
$D$-comodule structure. Hence we have

\begin{theorem}\label{eq:equiv-bundles}
If $C$ is a coalgebra and ${\cal E}$ is left $D$- right $G$-bicomodule,
then the equivalences in {\rm Theorems} 
\ref{thm:invarGlobSec}, \ref{thm:globalsec}
respect the $D$-comodule structure.
In particular, if ${\cal G} = D = {\cal E}$ is a Hopf algebra,
and ${\cal B}$ a quantum subgroup, then 
$\Gamma L_\chi$ is isomorphic
to the induced ${\cal G}$-comodule from ${\cal B}$.
\end{theorem}

\section{Background: Perelomov coherent states}
\label{sec:perelRawn}

Perelomov coherent states generalize 
the Schr\"{o}dinger coherent states 
to the Lie group setting~(\cite{perelomov, Perelomov:1972bd}).

We use the geometric language of~\cite{Rawnsley77}; cf. also~\cite{Ber74}.

Let $G^{\DDl C}$ be a complex connected semisimple Lie group with 
compact real form $G$, and a Borel subgroup $B$. 
We will often view these groups as affine algebraic groups over $\DDl C$.
Let $\chi : B \rightarrow {\DDl C}$ be a character of $B$
and ${\DDl C}_\chi$ the corresponding 1-dimensional $B$-module.
The projection $p : G^{\DDl C}\rightarrow G^{\DDl C}/B$ 
defines a principal $B$-bundle. The associated bundle
$L_\chi = G^{\DDl C} \times_\chi {\DDl C}_\chi$ 
with the projection $p_L : L_\chi \to G^{\DDl C}/B$. 
The left action of $G$ on $G^{\DDl C}$ induces an action of $G$ on $L_\chi$
and the formula $(g_* s)(x) = gs(g^{-1}x)$ 
defines an action of $G^{\DDl C}$ on the space 
$V_\chi = \Gamma L_\chi$ of holomorphic sections of $L_\chi$ 
which is by {\sc Borel-Weil} theorem, an irreducible unitarizable $G$-module. 
An invariant unitary product on $\Gamma L_\chi$,
antilinear in $1^{\rm st}$ and linear in $2^{\rm nd}$ argument, 
is denoted $\langle | \rangle$.
 
  Consider a (holomorphic) section $s \in \Gamma L_\chi$
and a nonzero point $q$ in some fiber ${p_L}^{-1}(x)$. Then
\[ s(x) = s(p_L(q)) = l_q(s) q, \]
for some number $l_q(s)$. The correspondence
\[\begin{array}{lcr}
s \mapsto l_q(s),  & \mbox{ } &
l_q : \Gamma L_\chi \rightarrow {\DDl C},
\end{array}\]
is a continuous linear functional.
Using Riesz's theorem, we infer the existence of an element
\[ e_q \in \Gamma L_\chi \mbox{ such that } 
\,\,l_q(s) = \langle e_q | s \rangle. \]
The vectors (sections) of the form $e_q \in \Gamma L_\chi$ are called
{\bf coherent vectors}. Corresponding projective classes are called
coherent states.

\begin{prop}
(\cite{Rawnsley77, Skoda:thes}) \label{prop:rawn}
(i) $e_{gq} = g_* e_q$ for all $g \in G^{\DDl C}$.

(ii) $e_{cq} = \bar{c}^{-1} e_q$ for all $c \in {\DDl C}$.

(iii) Coherent states i.e. the projective classes
of all coherent vectors belong to the same projective orbit.

(iv)  The set of all $e_q$ where $q \in {(p_L)}^{-1}(1_G B)$ 
agrees with the set (ray) of all
heighest weight vectors in $V_\chi$ for fixed $B$.

(v) The set of all $e_q$ where $q \in {(p_L)}^{-1}(u)$ for fixed 
$u \in G^{\DDl C}/B$ is the heighest weight space for some subgroup of
$G^{\DDl C}$ conjugated to $B$.
\end{prop}

\newtheorem{cor}{Corollary}
\begin{cor}\label{cor:cstateslocal} 
Let $U \subset G^{\DDl C}/B$ be an open set, $q \in L_\chi$ given, 
and $t : U \rightarrow G^{\DDl C}$ a section of
the principal $B$-bundle $p^{-1}(U)\to U$.
For each $g \in p^{-1}(U)$ there is a unique decomposition
$g = t(gB) b$ such that $b \in B$. 
The following ``homogeneity'' formula holds:
\begin{equation} \label{eq:cs_homogeneity}
 g e_q = \chi^{-1}(b) e_{t(gB)q}
\end{equation}
\end{cor}

\begin{pf} We have $g_*e_q =  t(gB)_* b_* e_q = t(gB)_*\chi^{-1}(b) e_q$
by~(\ref{eq:coh_state_2}); taking into account that 
$\chi^{-1}(b)$ is a scalar, this equals to
$\chi^{-1}(b)t(gB)_* e_q$, 
hence by {\it(i)} of the Proposition~\ref{prop:rawn},
also to $\chi^{-1}(b) e_{t(gB)q}$.
\end{pf}

\begin{defn}
The {\bf local family of coherent vectors} 
corresponding to the triple $(U,t,q)$ is the map
\begin{equation}\label{eq:localfamilycs} 
\begin{array}{lcr}
C_{(U,t,q)} : U \rightarrow V_\chi \equiv \Gamma L_\chi, &\mbox{ } &
 C_{(U,t,q)} : [g] \mapsto e_{t([g])q}. 
\end{array}
\end{equation}
\end{defn}
For any $w$ in the Weyl group $W$ of $G$, there is a Zariski open
subset $G^{\DDl C}_w \subset G^{\DDl C}$ consisting of all
$g \in G^{\DDl C}$ for which there exists 
(automatically unique) $w$-Gauss decomposition $g = wyb$ 
where $y \in G^{\DDl C}$ belongs to the unipotent subgroup 
of the opposite Borel $B'$, and $b \in B$. 
Set $G^{\DDl C}_w$ is also $B$-invariant, 
hence a total space of the restricted fibration over
a Zariski open subset $G^{\DDl C}_w/B \subset G^{\DDl C}/B$. 
Define the local section
$t_w : G^{\DDl C}_w/B \rightarrow G^{\DDl C}_w \subset G^{\DDl C}$
by $t_w([g]) = wy$ where $g = wyb$ as above.
We denote 
\[ C_w := C_{(w,v_0)} := C_{(G^{\DDl C}_w/B,t_w,q)}\]
where $v_0 = e_q$ is a fixed highest weight vector in $V_\chi$. 
The collection of maps $\{C_w, w \in W\}$ 
will be generalized to the quantum group setting below.
They can be viewed as $C_w \in {\cal O}(G^{\DDl C}_w/B) \otimes V_\chi$
where ${\cal O}(G^{\DDl C}_w/B)$ is the complex algebra of
all algebraic functions on $G^{\DDl C}_w/B$.

In this particular case, Corollary~\ref{cor:cstateslocal} becomes
\begin{prop}
If $g = wyb$ is the Gauss decomposition in $G_w$ then for all $g \in G$
\begin{equation}\label{eq:coh_state_2}
  gv_0 = \chi^{-1}(b) C_w(gB),
\end{equation}
and $C_w$ is the unique element in
${\cal O}(G^{\DDl C}_w/B) \otimes V_\chi$ for which this holds.
\end{prop}

\section{Quantum coherent states and localizations}
\label{sec:QCSloc}

\begin{defn} Let $\chi$ be a group--like element 
in a Hopf algebra ${\cal B}$,
and $(V,\rho)$ a right ${\cal B}$-comodule.
A {\bf $\chi$-coinvariant} in $V$ is
an element $v_\chi \in V$ such that $\rho v_\chi = v_\chi \otimes \chi$.
\end{defn}

Let $\pi : {\cal G}\rightarrow {\cal B}$ 
be a surjective homomorphism of Hopf algebras.
We say that $(\pi, {\cal B})$ is a {\bf quantum subgroup} of ${\cal G}$.
Every ${\cal G}$-comodule (resp. comodule algebra) $(V,\rho)$ 
is a ${\cal B}$-comodule (comodule algebra)
via $\rho_{\cal B} = ({\rm id} \otimes \pi) \circ \rho$.
In particular, ${\cal B}$ coacts on ${\cal G}$ by 
$({\rm id} \otimes \pi) \circ \Delta_{\cal G}$ and this coaction makes
${\cal G}$ a left ${\cal B}$--comodule algebra and similarily 
for the right coactions. In particular, ${\cal G}$ can be
viewed as left-right ${\cal B}$-${\cal G}$--bicomodule 
${}^{\cal B}{\cal G}^{\cal G}$. Hence to each ${\cal B}$-comodule
$V$ one can attach an {\bf induced ${\cal G}$--comodule} by the formula
${\rm Ind}^{\cal G}_{\cal B} V = V \Box^{\cal B} {\cal G}$. 
This defines the induction functor ${\rm Ind}^{\cal G}_{\cal B}$ 
which is left adjoint to the restriction functor 
$(V,\rho) \mapsto (V,\rho_{\cal B})$ (Frobenius reciprocity for comodules).

\begin{defn}
Let $(V,\rho)$ be any ${\cal G}$-comodule 
and $(\pi,{\cal B})$ a quantum subgroup of ${\cal G}$.
A {\bf weight covector} of weight $\chi$, 
is any $\chi$-coinvariant for ${\cal B}$--coaction
$\rho_{\cal B}$ in $V$ i.e.
\[ ({\rm id} \otimes \pi) \rho v_\chi = v_\chi \otimes \chi. \]
Let $(V_\chi,\rho) = {\rm Ind}_{\cal B}^{\cal G}{\DDl C}_\chi$ 
be the induced right ${\cal G}$-comodule induced from
the 1-dimensional comodule $z \mapsto z \otimes \chi$.
\end{defn}
\begin{defn} \label{def:bialg-realform}
A $*$-involution on a ${\DDl C}$-bialgebra $H$ is an antilinear 
map $* : H \to H$, for which $(ab)^* = b^* a^*$, 
$\Delta(a^*) = \sum a_{(1)}^* \otimes a_{(1)}^*$
and $\epsilon(a^*) = \overline{\epsilon(a)}$.
A pair $(H,*)$ is called a {\bf real form} of $H$.
\end{defn} \nodo{comment on antipode}

\begin{lem} 
(Schur's lemma for comodules) Let $C$ be a coalgebra over
${\DDl C}$ and $(V, \rho)$ a right $C$-comodule. If $(V,\rho)$ is
finite-dimensional and simple (no coinvariant subspaces),
then every $C$-comodule map $A : V \rightarrow V$ 
equals $\alpha \cdot{\rm id}_V$
for some $\alpha = \alpha_A \in {\DDl C}$.
\end{lem}

\begin{defn}~(\cite{KlimykSchmud})
Let $H$ be a Hopf $*$-algebra. 
An inner product $\langle \cdot  | \cdot \rangle$ on a
right $H$-comodule $V$ is a coinvariant inner product iff
\[ \langle w | z \rangle 1_H = 
        \sum \langle w_{(0)} | z_{(0)} \rangle  z_{(1)} w^*_{(1)} \]
An $H$-comodule which is a Hilbert space via a coinvariant inner product
will be called a right {\bf unitary $H$-comodule}.
\end{defn}

Consider a real form of a Hopf algebra ${\cal G}$ with the following data:

$\bullet$ (D1)
A surjective map of Hopf algebras $\pi : {\cal G}\rightarrow {\cal B}$.

$\bullet$ (D2)
 A group--like element $\chi \in {\cal B}$.

$\bullet$ (D3)
A coinvariant inner product on $V_\chi$.

$\bullet$ (D4)
A weight covector $v_\chi \in V_\chi$ with norm $1$.

$\bullet$ (D5)
 A Zariski local trivialization 
$\Lambda = \{ \lambda = (\iota_\lambda, S^{-1}_\lambda {\cal G}, 
\gamma_\lambda ) \}_{\lambda \in \Lambda}$
of ${\cal G}$ as a right ${\cal B}$-comodule algebra.

From now on, $V$ will be a comodule over the real form of ${\cal G}$
with fixed unitary equivalence $V\cong V_\chi$ which we often treat
as an identification.
Denote by  $V^{\rm triv}$ the trivial ${\cal G}$-comodule 
with the same underlying vector space as $V$.
\begin{defn}
Let (D1-5) be given and $\lambda \in \Lambda$.
A (Zariski-) {\bf local family of coherent vectors} in $\lambda$
or a {\bf polynomial coherent vector}\footnote{
Terminological remark. '{\it Polynomial}' because it is a polynomial 
in the generators of algebra ${\cal E}^{{\rm co}{\cal B}}$
of localized coinvariants 
decorated (tensored) with coefficients in Hilbert space. 
This terminology is 
occasionally used in physics literature (in the group case, as well as in
the quantum group examples, e.g.~\cite{HPT}, p.~1382).
}
in $\lambda$
is an element $C_\lambda \in V \otimes {\cal G}^{{\rm co}{\cal B}}_\lambda$ 
such that
\begin{equation}\label{eq:cohstatesDefn}
\rho_\lambda v_\chi = C_\lambda \gamma_\lambda (\chi) 
\end{equation}
holds in $V \otimes {\cal G}_\lambda$
where $\gamma_\lambda (\chi)$ on the right multiplies the second tensor
factor in $C_\lambda$ and $\rho_\lambda$ is the localized 
${\cal B}$-coaction $(\id\otimes\iota_\lambda)\rho$.
A {\bf global family of coherent vectors} is an 
element $C$ of 
$\Gamma(V^{\rm triv} \otimes L_\chi)$ such that
${\cal K}(C) = {\cal K}_\Lambda(C_{\Lambda}) = \rho v_\chi$ 
(for one, hence any, choice of $\Lambda$). 
Then $\kappa_\lambda(C) = \rho_\lambda v_\chi$.
\end{defn}
{\it Remark.} Equality~(\ref{eq:cohstatesDefn}) is a generalization of
the identity~(\ref{eq:coh_state_2}) and related to
Proposition 5.11 in~\cite{JurSt:coh}.
\begin{prop} The following are equivalent:

a) There exists a global family of coherent states $C$;

b) There exists a local trivialization $\Lambda$ of ${\cal E}$
such that a local family of coherent states $C_\lambda$
exists for each $\lambda$ in $\Lambda$;

c) For each local trivialization $\Lambda$ of ${\cal E}$
and each $\lambda$ in $\Lambda$ there exists a local family $C_\lambda$
of coherent states in $\lambda$.

Since ${\cal K}_\Lambda$ is a natural equivalence, if (a-c) are
true, then the global family is unique. The same for the local family 
in any {\it given} local trivialization.
\end{prop}
\begin{pf}
An exercise to the reader:
use the globalization lemma and the explicit description of ${\cal K}_\Lambda$.
Notice though that given only one $C_\lambda$ does not always suffice. Indeed, 
$C_\lambda \gamma_\lambda(\chi) \gamma_{\lambda'}(S\chi)$ is a candidate
for $C_{\lambda'}$, but it does not need to extend to 
an element in $V\otimes{\cal G}_{\lambda'}$ in general.
\end{pf}

{\it Let us extend the product
$\langle|\rangle$ on $V \cong V \otimes {\DDl C}
\subset V \otimes {\cal G}$ 
to a sesquilinear form }
\[ \langle|\rangle : (V \otimes {\cal G})\otimes V 
\rightarrow {\cal G},\,\,\,\,
\left\langle\left. \sum_o w_i \otimes g_i \right| v \right\rangle
:= \sum_i \langle w_i | v\rangle g_i,\]
and analogously define $\langle|\rangle_{\lambda\mu\ldots}$ 
(often skipping the subscripts) on $V \otimes {\cal G}_{\lambda\mu\ldots}$.
In particular, for any $v \in V_\chi$ the expression 
$\langle C_\lambda | v\rangle_\lambda$ is an element in ${\cal G}_\lambda$.
Let $|C_\lambda\rangle := C_\lambda$ in such context.

\begin{prop}
For each $v \in V$,
$\prod_\lambda \langle C_\lambda | v\rangle_\lambda$ 
is an element in $\Gamma_\Lambda L_\chi$,
and hence, by  Theorem~\ref{thm:globalsec}, it determines an
element in $V_\chi\cong V$.
\end{prop}
\begin{pf} 
By the definition, 
$\langle C_\lambda|v\rangle \in {\cal G}_\lambda^{{\rm co}{\cal B}}$.
Hence, by~(\ref{eq:gammalchi}), 
for each pair $(\lambda, \lambda')$, we have to check that 
$\langle C_\lambda|v\rangle_\lambda \gamma_\lambda (\chi)
= \langle C_{\lambda'}| v\rangle_{\lambda'} \gamma_{\lambda'}(\chi)$ 
in both consecutive localizations. 
To that aim observe that
\begin{equation}\label{eq:proofforfunctionals}
\langle\rho_\lambda v_\chi | v \rangle_\lambda = 
\langle C_\lambda \cdot (1 \otimes \gamma_\lambda \chi)
\,| v\rangle_\lambda = \langle C_\lambda | v\rangle_\lambda
\gamma_\lambda (\chi).
\end{equation}
Then observe, that symbol $\langle|v\rangle 
= \langle,v\rangle \otimes {\rm id}$ 
commutes with the localizations, in the sense that
$({\rm id} \otimes \iota^\lambda_{\lambda, \lambda'}) 
\circ \langle|v\rangle_\lambda
= \langle|v\rangle_{\lambda\lambda'} \circ 
({\rm id} \otimes \iota^\lambda_{\lambda, \lambda'})$.
Hence the equality 
$\rho_\lambda v_\chi = \rho_{\lambda'} v_\chi$, 
which may be fully expanded as
\begin{equation}\label{eq:rholam}
({\rm id} \otimes \iota^\lambda_{\lambda\lambda'}) 
({\rm id} \otimes \iota_{\lambda}) \rho v_\chi = 
({\rm id} \otimes \iota^{\lambda'}_{\lambda\lambda'})
({\rm id} \otimes \iota_{\lambda'})
\rho v_\chi
\end{equation}
implies that 
$\langle\rho_\lambda v_\chi | v\rangle 
= \langle\rho_{\lambda'} v_\chi | v\rangle$,
and by~(\ref{eq:proofforfunctionals}) this yields the wanted equality. 
The same way, using $\iota^{\lambda}_{\lambda'\lambda}$ and 
$\iota^{\lambda'}_{\lambda'\lambda}$ in~(\ref{eq:rholam}) this time,
check the identity in another consecutive localization.
\end{pf}

\section{Resolution of unity.}

A measure $\mu$ on a locally compact group $G$ is 
{\bf left-invariant} if
\[ \int_G f(gg') d\mu(g') = \int_G f(g') d\mu(g') \]
for all integrable functions $f$ on $G$ and for all $g \in G$.
Here we may replace $f(gg')$ by $(\Delta f) (g \otimes g')$ where $\Delta$
is the comultiplication in a suitable topological Hopf algebra of
functions on $G$. Evaluation at $g$ is a certain linear functional
$h_g$ on that algebra. This motivates~(\cite{Majid,Montg})
\begin{defn}
A {\bf left-invariant integral} (= left Haar integral) on a Hopf algebra $H$
is a linear functional $\int$ on $H$ such that
\[ \langle h \otimes \int, \Delta (f) \rangle = \langle  h, 1 \rangle
\langle \int, f \rangle , \,\,\,\,\,\,\,\,\,\, \forall h \in H^*.\]
A left Haar integral $\int$ is {\bf normalized} 
if $\langle \int , 1 \rangle = 1$.
\end{defn}
Since linear functionals separate elements of $H$, the left invariance
can be expressed as (dropping the evaluation brackets)
\[ ({\rm id} \otimes \int) \Delta(a) = ( \int a ) \cdot 1_H,
\,\,\,\,\,\,\,\,\,\,\,\forall a \in H.\]
In other words,
 $1_H \int$ intertwines $H$ (as a right $H$-comodule with respect to
the comultiplication) and its trivial subcomodule ${\DDl C} \cdot 1_H$.

\begin{theorem}\label{th:res-of-unity-comod}
Let $\int$ be a left integral on a Hopf $*$-algebra $H$, 
and $(V,\rho,\langle,\rangle)$
a simple unitary right $H$-comodule. Fix a vector $w \in V$.
Define the operator $A : V \rightarrow V$ by
\[ A|v\rangle 
=  \sum \langle w_{(0)} | v \rangle w_{(0)'} \int w_{(1)}^* w_{(1)'}\]
Then $A$ is a scalar operator. 
\end{theorem}

\begin{pf} In the following, the primed Sweedler indices belong to
another copy of the same variable, as in~\cite{Majid}. 
We compute directly
\[ \rho Av = \sum \langle w_{(0)} | 
v \rangle w_{(0)'} \int w_{(1)}^* w_{(2)'}\otimes w_{(1)'} \]
On the other hand,
\[ (A \otimes {\rm id}) \rho v = \sum
\langle w_{(0)} | v_{(0)} \rangle w_{(0)'} \int w_{(1)}^* w_{(1)'}
\otimes v_{(1)},\]
what is by the left invariance of the integral equal to
\[ \sum \langle w_{(0)} | v_{(0)} \rangle w_{(0)'} \int w_{(2)}^* w_{(2)'}
\otimes v_{(1)}w_{(1)}^* w_{(1)'},
\]
and, by the coinvariance of the inner product,
\[ (A \otimes {\rm id}) \rho v = \sum
\langle w_{(0)} | v \rangle w_{(0)'} \int w_{(2)}^* w_{(2)'}
\otimes w_{(1)'}\]

We conclude that $\rho Av = (A \otimes {\rm id})\rho v$. 
Hence the theorem follows from the Schur's lemma for comodules.
\end{pf}

\begin{defn}
$d\mu_\lambda(\chi) :=  \gamma_\lambda(\chi)(\gamma_\lambda(\chi))^*$ 
in ${\cal E}_\lambda$. 
\end{defn}

\begin{theorem}\label{eq:maintheorem}
Elements  $|C_\lambda \rangle d\mu_\lambda(\chi) \langle C_\lambda |
:= C_\lambda d\mu_\lambda(\chi) C^*_\lambda$
do not depend on $\lambda$ (agrees in all consecutive 
localization overlaps). Hence, by the globalization lemma,
this family defines localized representatives 
of a unique expression $|C\rangle d\mu (\chi) \langle C|$
in $V\otimes{\cal G}\otimes V^*$.
Taking a Haar integral in the tensor factor ${\cal G}$
yields a scalar operator $\alpha\cdot\id$
on $V$ (we identify ``states'' in
$V \otimes V^*$ with operators).
\end{theorem}

Remark: While $\langle C_\lambda | = C^*_\lambda$ should
live in ${\cal G}_\lambda\otimes V^*$, we may define it only as
a part of the expressions of the form $f^* C^*_\lambda := (C_\lambda f)^*$
with $f\in {\cal G}$ such that 
$C_\lambda f\in V\otimes \iota_\lambda({\cal G})$. 
Indeed, the involution $*$ is not defined on entire ${\cal G}_\lambda$,
but only on ${\cal G}$, or if you wish, $\iota_\lambda({\cal G})$.

\begin{pf} Notice that in each local trivialization $\lambda$,
\[
| C_\lambda \rangle\gamma_\lambda(\chi)(\gamma_\lambda(\chi))^* 
\langle C_\lambda | v \rangle =
\sum \langle w_{(0)} | v \rangle w_{(0)'} w_{(1)}^* w_{(1)'}
\]
where on the LHS we assume that the pairing between $V^*$ and $V$ is
assumed (applied) and on the RHS we assume appropriate localization. 
Recall that the product 
$C_\lambda \gamma_\lambda(\chi)$ does NOT
depend on the localization.  
By~(\ref{eq:gammalchi}) and 
Theorem~\ref{thm:globalsec} the RHS reads (no localizations
this time) the element in $V \otimes {\cal G}$ to integrate.
Hence by
Theorem~\ref{th:res-of-unity-comod} we see that
\[  \int C_\lambda  d\mu_{\lambda}(\chi) C^*_\lambda  = 
\int  | C_\lambda \rangle \gamma_\lambda(\chi)
(\gamma_\lambda(\chi))^*\langle C_\lambda |
\]
is a scalar operator.
\end{pf}

In other words, if a family $\{C_\lambda\}_\lambda$ 
of coherent states exists, then the
{\bf coherent states make a resolution of unity}. 
This fact enables us to define an analogue of 
the Bargmann transform~(\cite{Bargmann:I}).
To a vector $v \in V$ ($V$ is physically a space
describing some quantum numbers of the system; or a sector in
a decomposition of such a space) we assign 
$\langle C_\lambda | v \rangle \in {\cal G}_\lambda^{{\rm co}{\cal B}}$.
If $H$ is a linear operator on $V$, 
denote $H\,| C_\lambda \rangle:= (\id \otimes H) |C_\lambda \rangle$. 
Suppose $\alpha\neq 0$ is the constant from Theorem~\ref{eq:maintheorem}. 
Then
\[ H |v\rangle =
\alpha^{-1}\int  \,H\, | C_\lambda \rangle \,d\mu_{\lambda}(\chi)\,
\langle \,C_\lambda | v \rangle.\]
We then obtain (as in commutative case) a noncommutative version
of a reproducing ``integral'' kernel on a Hilbert space, and 
equations involving $H$ (e.g. deformations of Schr\"{o}dinger equation 
where $H$ is a Hamiltonian) can be
written down in this coherent state representation.

\section{Comments on the quantum group case}

The simplest example concerns the coherent states for
${\cal G} = {\cal O}(SU_q(2))$. We will mainly follow the notation
and conventions of~\cite{KlimykSchmud}. 
${\cal O}(SL_q(2))$ is a noncommutative Hopf algebra over ${\DDl C}$ 
with 4 generators $a,b,c,d$, usually assembled in a matrix 
$T = \left(\begin{array}{cc}a&b\\c&d\end{array}\right)$,
with relations $ab = qba$, $ac = qca$, $bc=cb$,
$bd = qdb$, $cd = qdc$, $ad-da = (q-q^{-1})bc$ 
and ${\rm det_q}T := ad - qbc = 1$.
${\cal O}(SU_q(2))$ is a real form of 
${\cal O}(SL_q(2))$ determined by formulas 
$a^* = d$, $b^* = -qc$, $c^* = -q^{-1}b$, $d^* = a$.
A vector space basis of ${\cal O}(SL_q(2))$ is
$\{ a^k b^r c^s\}_{k>0, r,s\geq 0} \cup 
\{ b^r c^s d^t\}_{r,s,t\geq 0}$. 
In particular, ${\cal O}(SL_q(2))$ splits into
a direct sum ${\DDl C}[\zeta]\oplus {\rm compl}(\zeta)$
where ${\DDl C}[\zeta]$ is the span of the basis elements of
the form $(bc)^r$ and ${\rm compl}(\zeta)$ the span of the rest of basis.
Notation ${\DDl C}[\zeta]$ suggests that it is the
algebra of polynomials in $\zeta = -qbc$, which will play
major role below. 
${\cal O}(SU_q(2))$ posses a unique Haar functional 
$\int$, found by {\sc Woronowicz}. With respect to the
direct sum decomposition above, $\int$
is nontrivial only on ${\DDl C}[\zeta]$ where it is given
by formulas involving {\sc Jackson}'s $q$-integral,
or equivalently~(\cite{KlimykSchmud}) 
$$\int \zeta^r = \frac{1 - q^{-2}}{1-q^{-2(r+1)}}, 
\,\,\,\,\,\,\,r = 0,1,2,\ldots$$
The lower quantum Borel subgroup ${\cal B}$ will be the quotient
${\cal O}(SL_q(2))/I$, 
where $I$ is the 2-sided ideal generated by $b$. 
$I$ is a Hopf ideal, hence ${\cal B}$ is a Hopf algebra. 
The quotient map $\pi : {\cal G}\to {\cal B}$ is datum (D1) 
from Sec.~\ref{sec:QCSloc}. The images of generators
are denoted $\lambda = \pi(a)$, $\xi = \pi(c)$,
$\lambda^{-1} = \pi(d)$ and $\pi(b) = 0$. 
Manin plane ${\cal O}({\DDl C}^2_q)$ is an algebra 
with two generators $x,y$ and a single
relation $xy = qyx$. Elements of the form $x^r y^s$ form
a basis of ${\cal O}({\DDl C}^2_q)$.
The latter is a right ${\cal O}(SL_q(2))$-comodule algebra via
$$\rho(x^r y^s) = (x\otimes a + y \otimes c)^r(x\otimes b + y\otimes d)^s.$$
${\cal O}({\DDl C}^2_q)$ splits into the homogeneous
components $V_n = \oplus_{r+s = n} {\DDl C} x^r y^s$ of dimension $n+1$,
which are irreducible and unitary.
Our datum (D2) will be $\chi = \lambda^{-n}$ in ${\cal B}$, 
(D3) $V_\chi = V_n$, and (D4) will be the weight vector $v_\chi = y^n$.
Datum (D5) is given by 1) two localizations ${\cal G}_b = {\cal G}[b^{-1}]$ and
${\cal G}_d = {\cal G}[d^{-1}]$ 
at Ore sets multiplicatively generated by $b$ and
$d$ respectively; 2) comodule algebra maps $\gamma_b,\gamma_d$ 
obtained from the quantum Gauss decomposition.
Let $u := bd^{-1}\in {\cal G}_d$. 
It is easy to show~(\cite{Skoda:thes, Skoda:loc-coinv1-and2}) 
that these localizations cover ${\cal G}$.
Both localizations are
$\rho_{\cal B}$-compatible, namely $\rho_{\cal B}$
extends by $\rho_{\cal B}(b^{-1}) = b^{-1}\otimes\lambda^{-1}$
and $\rho_{\cal B}(d^{-1}) = d^{-1}\otimes\lambda^{-1}$.
The algebras of localized $\rho_{\cal B}$-coinvariants
are given by ${\cal G}_b^{{\rm co}{\cal B}} = {\DDl C}[u]$
and ${\cal G}_d^{{\rm co}{\cal B}} = {\DDl C}[u']$
where $u = bd^{-1}$ and $u' = db^{-1}$. 
A unique (``Gauss'') decomposition of matrix $T$ in the form $wUA$ where 
$w$ is a permutation matrix, $U$ upper triangular unidiagonal and
$A$ lower triangular is possible in ${\cal G}_b$ with $w =\id$ 
and in ${\cal G}_d$ with 
$w = \left(\begin{array}{cc} 1 & 0 \\ 0 & 1 \end{array}\right)$.
Map $\lambda\mapsto A^1_1$, $\lambda^{-1}\mapsto A^2_2$
$\xi \mapsto A^2_1$ uniquely extend to a ${\cal B}$-comodule algebra 
map $\gamma_d : {\cal B}\to{\cal G}_d$, or to
$\gamma_b : {\cal B}\to{\cal G}_b$ in the latter case.
Explicitly $\gamma_d(\lambda) = a - bd^{-1}c$, 
$\gamma_d(\xi) = c$, $\gamma_d(\lambda^{-1}) = d$,
$\gamma_b(\chi_b) = d^n$; $\gamma_b(\xi) = c-db^{-1}a$,
$\gamma_b(\lambda) = a$, $\gamma_b(\lambda^{-1}) = b$.
Also $U^1_2 = u$ in ${\cal G}_d$ and $U^1_2 = u'$ in ${\cal G}_b$.

Analogously for general $n$, a cover of ${\cal O}(SL_q(n))$
by $n!$ $\rho_{{\cal B}_n}$-compatible Ore localizations $S_w$ 
($w$ in permutation group $\Sigma_n$) and ${\cal B}_n$-comodule 
algebra maps $\gamma_w : {\cal B}_n\to {\cal O}(SL_q(n))[S_w^{-1}]$ 
is a highly nontrivial fact which we have shown elsewhere.
It may be used to obtain the ${\cal O}(SU_q(n))$-coherent states. 

Using the $q$-binomial theorem, one obtains (in $V_n\otimes {\cal G}_d$)
$$
\rho(y^n) = \sum_{i = 0}^n {n \brack i}_{q^{-2}} 
q^{-{i \choose 2}} x^i y^{n-i} \otimes u^i d^n
$$
Basis vectors $v^n_i = \sqrt{{n \brack i}_{q^{-2}}} x^i y^{n-i}$
are orthonormal. Thus
$$
C_d :=  \sum_{i = 0}^n  q^{-{i \choose 2}} \sqrt{{n \brack i}_{q^{-2}}} 
v^n_i \otimes u^i,
$$
satisfies~(\ref{eq:cohstatesDefn}). Similar formula defines $C_b$
and the rest of requirements hold for these data. Thus
$$
A = \int_{SU_q(2)} \sum_{i,j=0}^n 
\sqrt{{n \brack i}_{q^{-2}} {n \brack j}_{q^{-2}}}
q^{{i\choose 2}+{j\choose 2}}
v^n_i\otimes (v^n_j)^* \otimes  u^i  d^n (u^j d^n)^*  
$$
\begin{lem}
$$\int_{SU_q(2)} u^i  d^n (u^j d^n)^* = \left\lbrace 
\begin{array}{ll}
0,\,\,\,\,\,\,\,\,&i\neq j\\ %
 {n \brack i}_{q^{-2}}^{-1} q^n q^{2{i\choose 2}}[n+1]^{-1}_q,
\,\,\,\,\,\,\,&i=j
\end{array}
\right.$$
\end{lem}
\begin{pf}  $(u^j d^n)^* = q^{i\choose 2}(d^*)^{n-j}(b^*)^j = 
q^{i\choose 2}(-q)^j a^{n-j} c^j$. The identity
$$d^r a^r = (1 + q^{-1}bc)(1+ q^{-3}(bc)^2)\ldots(1+q^{-2n-1}(bc)^r)
= (q^{-2}\zeta;q^{-2})_r,$$ 
implies $d^{n-i} a^{n-j} = d^{j-i} (q^{-2}\zeta;q^{-2})_{n-j}$
for $j\geq i$ and $(q^{-2}\zeta;q^{-2})_{n-j}a^{i-j}$ for $i<j$.
Thus, for $j\geq i$,
$$\begin{array}{lcl} u^i  d^n (u^j d^n)^* 
&=& q^{{i\choose 2}+{j\choose 2}}(-q)^j b^i d^{j-i}(q^{-2}\zeta;q^2)_{n-j} c^j 
\\&=& q^{{i\choose 2}+{j\choose 2}}(-q)^jb^i(q^{-2}\zeta;q^2)_{n-j} 
d^{j-i},\end{array}$$ 
what is for $j> i$ an element in ${\rm compl}(\zeta)$
hence it vanishes after integration, likewise an expression for
$i < j$, and only the terms with $i=j$ survive. Then
$u^i  d^n (u^i d^n)^* = - q^{2{i\choose 2}}\zeta^i (q^{-2}\zeta;q^{-2})_{n-i}$,
and using (52') in Chapter 4 of~\cite{KlimykSchmud} one derives
\begin{equation}\label{eq:ourqbeta}
\int_{SU_q(2)} \zeta^i (q^{-2}\zeta;q^{-2})_{n-i} = 
{n \brack i}_{q^{-2}} q^n [n+1]^{-1}_q,\end{equation}
and the rest of the calculation is immediate.
\end{pf}

Now
$$
A = - \int_{SU_q(2)} \sum_{i=0}^n {n \brack i}_{q^{-2}}
| i\rangle\langle i | \otimes \zeta^i (q^{-2}\zeta;q^{-2})_{n-i} 
$$

$$
\int_{SU_q(2)}  \sum_{i=0}^n {n \brack i}_{q^{-2}}| i\rangle\langle i |
\otimes \zeta^i (q^{-2}\zeta;q^{-2})_{n-i} = 
[n+1]^{-1}_q q^{-n} \sum_{i=0}^n | i\rangle\langle i |.
$$
The sum on RHS is of course the unity. The fact that there was no
additional factors depending on $i$ is the nontrivial property
of coherent states (Theorem~\ref{eq:maintheorem}). 
There are many proposals for ``$SU_q(2)$-coherent 
states'' in literature (search e.g. MathSciNet)
with similar (partly guessed) formulas with wrong $q$-factors
and still having some ``resolution of unity'' formulas. 
The wrong factors are compensated by 
effectively changing the measure as well, 
for which there is no freedom as $SU_q(2)$ has only one invariant integral
up to an {\it overall} constant! 

In other words, $\alpha = q^n [n+1]^{-1}_q$ and
the resolution of unity is 
$$
I = q^{-n} [n+1]_q \int_{SU_q(2)} | C\rangle d\mu(\chi) \langle C|.
$$

Formula~(\ref{eq:ourqbeta}) boils down to an integral representation of
{\sc Ramanujan}'s $q$-beta function (Theorem 10.3.1 in~\cite{AAR:special};
cf. also~\cite{Askey})
$$
\int_0^1 x^{\alpha} \frac{(qx;q)_\infty}{(q^\beta x;q)_\infty} d_q x = 
\frac{\Gamma_q(\alpha)\Gamma_q(\beta)}{\Gamma_q(\alpha+\beta)}.
$$
If $\beta \geq 1$ is a positive integer, then
the ratio in the integrand equals a polynomial in $q$ and $x$, namely,
$(1-x)(1-qx)\ldots(1-q^{\beta-1}x)$. Like for the ordinary 
beta function, there is another integral representation
involving $q$-integral from $0$ to $\infty$ 
with a polynomial in the {\it denominator}.
Namely, instead of the Haar integral over $\zeta = -qbc$
one effectively has a geometric integration over (deformed)
2-sphere with real coordinates $u = bd^{-1}$ and $\bar{u}$. 
However, in the denominator form, 
new $q$-factors appear depending on $q^i$.  
{\sc Jur\v{c}o}~(\cite{Jurco:1991tg})
wrote a similar formula without extra $q^i$ factors,
but both the ``measure'' and the coherent states are changed.
Hence those ``coherent states'' do not satisfy the 
defining factorization property~(\ref{eq:cohstatesDefn}) 
and the measure is not the invariant one. 

Computations of coherent states in selected local coordinates
in concrete examples  ${\cal O}(SU_q(n))$ with $n = 2,3$,
appeared in~\cite{HPT, Jurco:1991tg, Sazdjian:1995yg}, though
without full geometric justification, and sometimes with nongeometric
factors. 
Rudiments of another picture involving quantum group coherent states, 
related to geometric  quantization and orbit method, 
are discussed in~\cite{Soib:orbit}.
Finally, a local picture (i.e. calculations in main Bruhat cell) of 
the coherent states for the case of compact forms 
of quantum groups of types A,B,C,D, 
which differs from but is related to ours,
is in impressive work by {\sc Jur\v{c}o} and 
{\sc \v{S}\v{t}ov\'{i}\v{c}ek}~(\cite{JurSt:dres, JurSt:coh}). 
Their family of coherent states, $\Gamma$ (cf. (5.1) in ~(\cite{JurSt:coh})),
live in $V\otimes {\cal G}$, i.e. generalize a map 
$G\to V$ rather than $G/B\to V$. They however calculate some 
expression in corresponding coordinates on a cell in homogeneous
space, working in a big Zariski open cell (without 
rigorous justification for localization). 
Proposition 5.11 in~(\cite{JurSt:coh}) is stating the factorization
property (our formula~(\ref{eq:cohstatesDefn}))
of their quantity $w_\lambda^{-1} \langle \Gamma, u\rangle$
which ``belongs to some completion'' and basically agrees with our
coherent states. Their construction relies on structure properties 
of quantum groups, while our axiomatics allows a priori
treatment of Hopf algebras of more general origin. 
Furthermore, our construction utilizes the globalization of the 
geometry on the quantum homogeneous space.

In commutative case, the elements of a family of coherent states form
the projective orbit of the highest weight vector. 
The generators of ${\cal G}_{\lambda}^{{\rm co}{\cal B}}$ 
are the analogues of the local coordinates on 
a big open cell in the coset space, and
the coherent vector $C_\lambda$ may be viewed as a parametrization
of an open set in projective orbit by points in a coset space.
In similar spirit, in the case of
${\cal O}(SL_q(3))$, the reference~\cite{Sazdjian:1995yg}
views ${\cal G}_\lambda^{{\rm co}{\cal B}}$
as an analogue of the (algebra of functions on) unipotent group
parametrizes quantum orbit 
(though they note this algebra is not a bialgebra, 
unlike the classical case).
Here we clarify that, as in the
classical case, this should be understood as a parametrization
of an open dense subset of orbit, 
the latter being a noncommutative space.

\section{Open question: minimal uncertainty}

It remains to study {\bf ``covariant minimal uncertainty''}
properties of these coherent states. In the classical case, 
there is a quantity $\Delta({\frak C})$ which is a sort of a
``dispersion'' of the Casimir element ${\frak C}$, 
and it is minimized on the coherent state orbit.
We can show a quantum version by direct computation 
in one very simple example,
but it remains to be studied in greater generality.
An {\it Ansatz} for expression $\Delta({\frak C}_q)$, 
which attains minimum at highest weight 
and in limit $q\to 1$ gives $\Delta({\frak C})$,
has been proposed for the case of the standard quantum groups
possessing $R$-matrix in~\cite{Delbourgo:as}. One may hope to
reinterpret their expression in terms of braided
Casimir element~\cite{gur:index}, 
and then place it into our geometrical context.
In noncommutative case, {\it individual} coherent states
are not defined, but the noncommutative 'family' of coherent states
as a whole is still defined.
Hence one may try to show that the average of the ``dispersion''
over the {\it family of coherent states}
agrees with the minimal value over the whole representation space. 
In commutative case, this property is equivalent to saying that
the minimum is achieved at ``almost every'' 
(in the measure sense, what forces ``every'' by continuity)
point of the coherent state orbit.

Let me remark on possible strategies
to obtain the minimization property.
In the commutative case, the property 
can be traced to symplectic geometry. 
Up to a constant shift, and a negative multiple,
the quantity $\Delta({\frak C})$ equals the value of ``the square
of moment map'' $\|\mu\|^2$, properly understood~(\cite{Spera}). It is
essential that the coherent state orbits are K\"{a}hler, hence symplectic.
The square of moment map is extremal at symplectic orbits, 
hence the minimization of $\Delta({\frak C}_q)$~(\cite{Spera}). 
The noncommutative symplectic geometry developed by
{\sc Kontsevich, Ginzburg} and others~(\cite{ginzburg:symp,Konts:formalsymp})
may suggest a path to extend the moment map argument to
the noncommutative case. 

A noncommutative infinitesimal neighborhood 
of a commutative scheme has been introduced
by {\sc Kapranov}~(\cite{Kapranov}). It is a very flexible
setup for geometry with touch with direct calculations
({\sc Feynman-Maslov} calculus).
In the case of $q$-deformed groups one may hope to combine his
approach with filtration arguments 
(with finite-dimensional geometry of graded pieces), 
and extend Morse-like arguments for study of extremality.

\begin{acknowledgements}
I thank Prof. {\sc Baha Balantekin} for suggesting the study of coherent
states for quantum groups in a semesteral project in Fall 1996, 
and for his encouragement to try finding the ``measure'',
what later resulted in $d\mu(\chi)$. 
I thank Prof. {\sc Joel W. Robbin} for his patience, critique and many 
discussions on orbit method, coherent states and quantum groups
at the University of Wisconsin-Madison.
I thank the organizers of the 1999. summer workshop ``Geometry and
Physics'' at Santa Barbara during which,
in peace of my ITP office, Theorem~\ref{th:res-of-unity-comod}
has been proved (and an imprecise 
version of Theorem~\ref{eq:maintheorem}) while supported
in part by NSF under grant No. PHY99-07949. The project has been
finished at the Indiana University, Bloomington, and the final version of the
paper at the institutes Rudjer Bo\v{s}kovi\'{c}, Zagreb, 
and Max Planck, Bonn.
\end{acknowledgements} 


\end{document}